\documentclass[12pt]{article}
\usepackage[cp1251]{inputenc}
\usepackage[english, russian]{babel}
\usepackage{a4wide}
\usepackage{amssymb}
\usepackage{amsfonts}
\usepackage{amsmath}

\begin{document}

\begin{center}
{\bf Cauchy Problem for a Linear System of Ordinary Differential Equations 
of  the Fractional Order}
\vskip 3mm
{\bf Murat Mamchuev}%
\vskip 2mm
\end{center}
{{\footnotesize \quad

{\bf Abstract.}
The paper investigates the initial problem for a linear system of ordinary differential equations with the fractional differentiation operator Dzhrbashyan -- Nersesyan with constant coefficients.
The existence and uniqueness theorems of the solution of the boundary value problem under study are proved. 
The solution is constructed explicitly in terms of the Mittag-Leffler function of the matrix argument. 
The Dzhrbashyan -- Nersesyan operator is a generalization of the Riemann -- Liouville, Caputo and Miller-Ross fractional differentiation operators. 
The obtained results  as special cases contain results related to the study of initial problems for systems of ordinary differential equations with Riemann -- Liouville, Caputo and Miller -- Ross derivativess, and the investigated initial problem generalizes them.

}}%

\vskip 7mm

\begin{center}
{\bf Introduction}
\end{center}

Consider the system of ordinary differential equations
\begin{equation} \label{sn}
Lu(x)\equiv D_{0x}^{\{\alpha_0, \alpha_1, ..., \alpha_m\}} u(x)-Au(x)=f(x),
\end{equation}
where $D_{0x}^{\{\alpha_0, \alpha_1, ..., \alpha_m\}}$  is the
Dzhrbashyan -- Nersesyan fractional differentiation operator of the order
$\alpha=\sum\limits_{i=0}^{m}\alpha_i-1>0$
\cite{Djrb-Ners-1968},
 $\alpha_i \in (0,1]$ $(i=\overline{0,m});$
$f(x)= ||f_1(x),f_2(x),...,f_n(x)||$  and
$u(x)= ||u_1(x),u_2(x),...,u_n(x)||$
are respectively the given and unknown $n$-vectors,
$A$ is a given constant  $n \times n$ matrix.

The Dzhrbashyan -- Nersesyan fractional differentiation operator $D_{st}^{\{\gamma_0, \gamma_1,...,\gamma_m\}}$
associated with the sequence $\{\gamma_0, \gamma_1,...,\gamma_m\},$ of order
$\gamma=\sum\limits_{i=0}^{m}\gamma_i-1>0,$ $\gamma_i\in(0,1],$ $(i=\overline{0,m})$
determined by the ratio \cite{Djrb-Ners-1968}
$$D_{st}^{\{\gamma_0, \gamma_1,...,\gamma_m\}}v(t)=D_{st}^{\gamma_m-1}D_{st}^{\gamma_{m-1}}...D_{st}^{\gamma_1}D_{st}^{\gamma_0}v(t),$$
$$D_{st}^{\{\gamma_0\}}v(t)=D_{st}^{\gamma_0-1}v(t),$$
where $D_{st}^{\nu}$ is the Riemann -- Liouville fractional integro-differentiation operator of order $\nu.$ 
The operator $D_{st}^{\nu}$ for $\nu< 0$ is defined as follows \cite [p. 9]{Nakhushev-2003}:
\[D_{st}^{\nu}v(t)=\frac{\mathop{\rm sgn}(t-s)}{\Gamma(-\nu)}
\int\limits_s^t \frac{v(\xi)d\xi}{|t-\xi|^{\nu+1}},\]
where $\Gamma(z)$ is the Euler gamma function.
For $\nu \geq 0 $ the operator $ D_{st}^{\nu} $ can be determined using the recursive relation
\[D_{st}^{\nu}v(t)=\mathop{\rm sgn}(t-s)\frac{d}{dt}D_{st}^{\nu-1}v(t).
\]

In 1954, J.H. Barrett \cite{Barrett-1954} investigated the initial problem for the equation
$$D_{ax}^{\alpha}y(x)+\lambda y(x)=h(x),$$
$$\lim\limits_{x\to a}D_{ax}^{\alpha-i}y(x)=K_i, \quad i=1,2,...,n,  \quad n-1<{\rm  Re}(\alpha)<n,   \quad n\in \mathbb{N}.$$

In 1968, M.M. Dzhrbashyan  and A.B. Nersesyan \cite{Djrb-Ners-1968} introduced the fractional differentiation operator
$D_{0t}^{\{\gamma_0, \gamma_1,...,\gamma_m\}}$ and
investigated the Cauchy problem for the equation (\ref{sn}) for $n = 1. $

Systems of linear ordinary equations of fractional order were first investigated in the works of V.K. Veber \cite{Veber-1976} -- \cite{Veber-1985} and
M.I. Imanaliev, V.K. Veber \cite{ImVeb-1980}.
In   1976, V.K. Veber \cite{Veber-1976}
in terms of the Mittag-Leffler function of the matrix argument, the solution of the Cauchy problem for the system of equations
$$D_{0x}^{\alpha}y(x)=Ay(x), \quad 0<\alpha\leq 2$$
with constant matrix $A$ was written.
The asymptotic behavior as $ x \to \infty $ of various solutions of this system (including the fundamental matrix) was studied in \cite{Veber-1983} and \cite{ImVeb-1980}.
Later V.K. Veber considered the Cauchy problem for an inhomogeneous system
$$D_{0x}^{\alpha}y(x)=A(x)y(x)+f(x),\quad n-1<\alpha\leq n,\quad n=1,2,...$$
with continuous matrix function $A(x)$ for $x\geq 0$ \cite{Veber-1985-1},
and in \cite{Veber-1985} he constructed a fundamental solution of this system with a constant matrix $A$ in terms of the Mittag-Leffler function of a matrix argument.
Also in the paper \cite{Veber-1985} are examples of some applications lead to systems of equations with fractional derivatives.

A.A. Chikriy and I.I. Matychyn in \cite{Matychyn-2007} and \cite{Matychyn-2008}, using the Laplace transform, obtained solutions of Cauchy problems for systems of equations of the form (\ref{sn}), with Riemann -- Liouville, Caputo and Miller -- Ross derivatives.

In the works \cite{Matychyn-2018-1} and \cite{Matychyn-2018-2} I. Matychyn and V. Onyshchenko, using the Mittag-Leffler matrix function, analytically and numerically investigated the solutions of the initial problems for systems of equations with fractional Riemann -- Liouville and Caputo derivatives.

Note that \cite{Mamchuev-2008-2} and \cite{Mamchuev-2019} investigated boundary value problems for multidimensional systems of partial differential equations of fractional order.
In \cite{Mamchuev-2013-m}, attention was drawn to the fact that in the one-dimensional case these results coincide with the results of \cite{Veber-1976}.

In this paper, we investigate the initial problem for the system (\ref{sn}).
We prove some properties of the Mittag-Leffler matrix function, obtain a general representation of the solutions to the system (\ref{sn}), and prove the theorem on the unique solvability of the Cauchy problem for this system.

\begin{center}
{\bf 1. Statement of the problem and the solvability theorem}
\end{center}


{\bf Problem 1.}
{\it Find a solution $u(x)$ of the system (\ref{sn}) in the interval $(0,l),$ with the conditions
\begin{equation} \label{eq110}
\lim\limits_{x \rightarrow 0} D_{0x}^{\{\alpha_0, \alpha_1,...,\alpha_k\}} u(x)=u_0^k,\quad
 0\leq k\leq m-1,
\end{equation}
where $u_0^k=||u_1^k, u_2^k,...,u_n^k||$ are given real constant $n$-vectors. }

A solution $u(x)$ of the system (\ref{sn}) such that
$D_{0x}^{\{\alpha_0, \alpha_1,...,\alpha_m\}}u(x)\in C(0,l),$
$D_{0x}^{\{\alpha_0, \alpha_1,...,\alpha_k\}}u(x)\in C^1(0,l)\cup C[0,l]$
$(k=\overline{0,m-1}),$  is called a regular solution of the system (\ref{sn}) in the interval $(0,l).$

{\bf Theorem 1.}
{\it Let $\alpha_0+\alpha_m>1,$ function $f(x)\in C(0,l)$ can be presented as $f(x)=D_{0x}^{\alpha_m-1}f_0(x),$
where $f_0(x)\in L(0,l).$
Then there exists a unique regular in the interval $(0,l)$ solution of Problem 1. Solution can be
represented as
\begin{equation} \label{eq111}
u(x)=\sum\limits_{i=1}^{m}D_{0x}^{\{\alpha_m, \alpha_{m-1}, ..., \alpha_{m+1-i}\}}G(x)u_0^{m-i}
+\int\limits_0^xG(x-t)f(t)dt,
\end{equation}
where 
$G(x)=x^{\alpha-1}E_{\alpha, \alpha}\left(Ax^{\alpha}\right).$
}

\begin{center}
{\bf 2. Preliminaries}
\end{center}


The formula
\begin{equation} \label{eq11}
\int\limits_{0}^{x} h(t)D_{0t}^{\nu}g(t)dt=\int\limits_{0}^{x} g(t)D_{xt}^{\nu}h(t)dt, \quad \nu<0
\end{equation}
is known as a formula of fractional integration by parts \cite [p. 9]{Nakhushev-2003}.

The following formula
$$D_{0x}^{\nu}\int\limits_{0}^{x} h(x-t)g(t)dt=\int\limits_{0}^{x} h(x-t)D_{0t}^{\nu}g(t)dt+$$
\begin{equation} \label{eq13}
+h(x)\lim\limits_{t\to 0}D_{0t}^{\nu-1}g(t), \quad 0<\nu<1
\end{equation}
holds \cite[p. 99]{Podlubny-1999}.

The following formula of fractional integro-differentiation of power function 
\begin{equation} \label{eq12}
D_{ax}^{\nu}\frac{|x-a|^{\mu-1}}{\Gamma(\mu)}=\frac{|x-a|^{\mu-\nu-1}}{\Gamma(\mu-\nu)}
\end{equation}
is valid for $\mu>0$ if $\nu\in {\mathbb R},$ and for $\mu \in {\mathbb R}$ if $\nu\in {\mathbb N}.$

In 1903, Mittag-Leffler introduced the function \cite{Mittag-Leffler-1903} $E_{\alpha}(z)$
\begin{equation} \label{eq141}
E_{\alpha}(z)=\sum\limits_{k=0}^{\infty}\frac{ z^k}{\Gamma(\alpha k+1)}, \quad
(\alpha \in {\mathbb C}; {\rm Re} (\alpha)>0),
\end{equation}
which is now known as the Mittag-Leffler function.

In 1905, A. Wiman \cite{Wiman-1905} generalized this function with the two-parameter Mittag-Leffler function
$E_{\alpha, \beta}(z)$
(also sometimes called the Wiman's function)
$$E_{\alpha, \beta}(z)=\sum\limits_{k=0}^{\infty}\frac{ z^k}{\Gamma(\alpha k+\beta)}, \quad
(\alpha, \beta \in {\mathbb C}; {\rm Re} (\alpha)>0,  {\rm Re} (\beta)>0).$$

The following properties of this function are valid:
\begin{equation} \label{eq15}
E_{\alpha,\beta}(z)=\frac{1}{\Gamma(\beta)}+zE_{\alpha,\beta+\alpha}(z),
\end{equation}
\begin{equation} \label{eq16}
D_{0x}^{\mu}x^{\beta-1}E_{\alpha,\beta}(\lambda x^{\alpha})=
x^{\beta-\mu-1}E_{\alpha,\beta-\mu}(\lambda x^{\alpha}), \quad (\beta>0, \, \mu\in {\mathbb R}).
\end{equation}

In 1971, an even more general function $E_{\alpha, \beta}^{\gamma}(z)$ was introduced by T.R. Prabhakar  \cite{Prabhakar-1971}   
$$E_{\alpha, \beta}^{\gamma}(z)=\sum\limits_{k=0}^{\infty}\frac{(\gamma)_k}{\Gamma(\alpha k+\beta)}\frac{ z^k}{k!}, \quad
(\alpha, \beta, \gamma \in {\mathbb C}; {\rm Re} (\alpha)>0,  {\rm Re} (\beta)>0, {\rm Re} (\gamma)>0 ),$$
where
$(\gamma)_k=\frac{\Gamma(\gamma+k)}{\Gamma(\gamma)}$
is the Pohhammer symbol. 
This function is known as the Prabhakar function.


The following equalities hold \cite{Prabhakar-1971}, \cite{Shukla-2007}:
\begin{equation} \label{eq180}
E_{\alpha,\beta}^1(z)=E_{\alpha,\beta}(z),
\end{equation}
\begin{equation} \label{eq18}
E_{\alpha,\beta}^{\gamma}(z)=E_{\alpha,\beta}^{\gamma-1}(z)+zE_{\alpha,\beta+\alpha}^{\gamma}(z),
\quad ({\rm Re} (\alpha)>0, \, \beta>0, \, \gamma\in {\mathbb C}),
\end{equation}
\begin{equation} \label{eq19}
D_{0x}^{\mu}x^{\beta-1}E_{\alpha,\beta}^{\gamma}(\lambda x^{\alpha})=
x^{\beta-\mu-1}E_{\alpha,\beta-\mu}(\lambda x^{\alpha}), \quad (\beta>0, \, \mu\in {\mathbb R}),
\end{equation}
\begin{equation} \label{eq155}
\lim\limits_{z \to 0}E_{\alpha,\beta}^{\gamma}(z)=\frac{1}{\Gamma(\beta)}.
\end{equation}

In 1976, V.K. Veber \cite{Veber-1976} introduced the Mittag-Leffler function of the matrix argument
(see also \cite{Matychyn-2018-1}).
Note that in the paper \cite{Veber-1976} naturally come into sight a function $E_{\alpha, \beta}^{\gamma}(z)$ with natural values of the parameter $ \gamma \in{\mathbb N}.$

Here we give a definition of the Mittag-Leffler function of a matrix argument, and then examine some of its properties.

Let $A$ be a square matrix and $H$ a matrix reducing the matrix $A$ to Jordan matrix i.e.
$$A=HJH^{-1}=H\big[J_{r_1}(\lambda_1), J_{r_2}(\lambda_2), ... , J_{r_p}(\lambda_p)\big]H^{-1},$$
where
$J(\lambda)={\rm diag}[J_1(\lambda_1),...,J_p(\lambda_p)]$ is a quasi-diagonal matrix with blocks of the form
\begin{equation*}
J_{k}\equiv J_{k}(\lambda_k)= \left \|
\begin{array}{cccl}
\lambda_k 	& 1   		& \ldots & 0  \\
    		& \lambda_k 	& \ldots & 0 \\
   		&  0  		& \ddots & \vdots \\
   		&     		&        & \lambda_k
\end{array}
\right \|,  \quad  k=1,...,p,
\end{equation*}
$\lambda_1,...,\lambda_p$ are the eigenvalues of the matrix $A,$
$J_k(\lambda_k)$ are square matrices of order $r_k+1,$
$\sum\limits_{k=1}^{p}r_k+p=n.$
Then we have
\begin{equation} \label{eq191}
E_{\alpha, \beta}(Az)=HE_{\alpha, \beta}(Jz)H^{-1}=H{\rm diag}\big[E_{\alpha, \beta}[J_1(\lambda_1)z], ... , E_{\alpha, \beta}[J_p(\lambda_p)]\big]H^{-1},
\end{equation}
where
\begin{equation*}
E_{\alpha, \beta}[J_k(\lambda)z]= \left \|
\begin{array}{cccl}
e_0 	& e_1   	& \ldots 	& e_{r_k-1} \\
    	& e_0 		& \ldots 	& e_{r_k-2} \\
		&  0  		& \ddots 	& \vdots \\
   		& 			&        	& e_0
\end{array}
\right \|,
\end{equation*}
$$e_n\equiv e_n^{\alpha,\beta}(\lambda,z)=z^n\sum\limits_{k=0}^{\infty}\frac{(k+n)!}{k!n!}\frac{\lambda^k z^k}{\Gamma(\alpha (k+n)+\beta)}=z^n E_{\alpha, \alpha n+\beta}^{n+1}(\lambda z).$$
It's obvious that
$$e_0\equiv e_0^{\alpha,\beta}(\lambda,z)=E_{\alpha, \beta}(\lambda z).$$

\begin{center}
{\bf 3. Properties of the Mittag-Leffler function of a matrix argument}
\end{center}

1. The following equality holds:
\begin{equation} \label{eq51}
E_{\alpha, \beta}(Az)=\frac{1}{\Gamma(\beta)}I+AzE_{\alpha, \beta+\alpha}(Az), \quad (\beta>0, \, \mu\in {\mathbb R}).
\end{equation}

Indeed, by the formula (\ref{eq15}), we obtain
$$e_0^{\alpha,\beta}(\lambda,z)=E_{\alpha,\beta}(\lambda z)=\frac{1}{\Gamma(\beta)}+\lambda zE_{\alpha,\beta+\alpha}(\lambda z)=
\frac{1}{\Gamma(\beta)}+\lambda z e_0^{\alpha,\beta+\alpha}(\lambda,z).$$
For $ n\not=0,$ taking into account the formula (\ref{eq18}), we have
$$e_n^{\alpha,\beta}(\lambda,z)=z^n E_{\alpha, \alpha n+\beta}^{n+1}(\lambda z)=
z^n E_{\alpha, \alpha n+\beta}^{n}(\lambda z)+\lambda z^{n+1} E_{\alpha, \alpha n+\beta+\alpha}^{n+1}(\lambda z)=$$
$$=ze_n^{\alpha,\beta+\alpha}(\lambda,z)+\lambda ze_{n+1}^{\alpha,\beta+\alpha}(\lambda,z).$$
From the last equalities we obtain
$$\left \|
\begin{array}{cccl}
e_0 	& e_1   	& \ldots 	& e_{r_k-1} \\
    	& e_0 		& \ldots 	& e_{r_k-2} \\
		&  0  		& \ddots 	& \vdots \\
   		& 			&        	& e_0
\end{array}
\right \|=\frac{1}{\Gamma(\beta)}I+z
\left \|
\begin{array}{cccl}
\lambda\bar{e}_0 		& \bar{e}_1+\lambda \bar{e}_1 	& \ldots 		& \bar{e}_{r_k-1}+\lambda \bar{e}_{r_k-2} \\
    	    			& \lambda\bar{e}_0 				& \ldots 		& \bar{e}_{r_k-2}+\lambda \bar{e}_{r_k-2}\\
						&  0  							& \ddots 		& \vdots \\
   						& 								&        		& \lambda \bar{e}_0
\end{array}
\right \|=$$
$$=\frac{1}{\Gamma(\beta)}I+z\left \|
\begin{array}{cccl}
\lambda_k 	& 1   		& \ldots & 0  \\
    		& \lambda_k 	& \ldots & 0 \\
   		&  0  		& \ddots & \vdots \\
   		&     		&        & \lambda_k
\end{array}
\right \|
\left \|
\begin{array}{cccl}
\bar e_0 	& \bar e_1   	& \ldots 	& \bar e_{r_k-1} \\
    		& \bar e_0 		& \ldots 	& \bar e_{r_k-2} \\
			&  0  			& \ddots 	& \vdots \\
   			& 				&        	& \bar e_0
\end{array}
\right \|,
$$
where $\bar e_n\equiv e_n^{\alpha,\beta+\alpha}(\lambda,z)=z^n E_{\alpha, \alpha n+\alpha+\beta}^{n+1}(\lambda z).$
Thus,
$$E_{\alpha, \beta}[J_{r_k}(\lambda_k)z]=\frac{1}{\Gamma(\beta)}I+J_{r_k}(\lambda_k) z E_{\alpha, \beta+\alpha}[J_{r_k}(\lambda_k)z].$$
Hence, by virtue of (\ref{eq191}), the equality (\ref{eq51}) follows.


2. For the matrix function, an analogue of the formula (\ref{eq16}) of fractional integro-differentiation is valid
\begin{equation} \label{eq52}
D_{0x}^{\mu}x^{\beta-1}E_{\alpha, \beta}(Ax^{\alpha})=x^{\beta-\mu-1}E_{\alpha, \beta-\mu}(Ax^{\alpha}), 
\quad (\beta>0, \, \mu\in {\mathbb R}).
\end{equation}

Indeed, by virtue of the equality (\ref{eq19}) we have
$$
D_{0x}^{\mu}x^{\beta-1}E_{\alpha, \beta}[J_{r_k}(\lambda_k)z]=
D_{0x}^{\mu}x^{\beta-1}\left \|
\begin{array}{cccl}
e_0 	& e_1   	& \ldots 	& e_{r_k-1} \\
    	& e_0 		& \ldots 	& e_{r_k-2} \\
		&  0  		& \ddots 	& \vdots \\
   		& 			&        	& e_0
\end{array}
\right \|=$$
$$=x^{\beta-\mu-1}
\left \|
\begin{array}{cccl}
\tilde e_0 	& \tilde e_1   	& \ldots 	& \tilde e_{r_k-1} \\
    		& \tilde e_0 		& \ldots 	& \tilde e_{r_k-2} \\
			&  0  			& \ddots 	& \vdots \\
   			& 				&        	& \tilde e_0
\end{array}
\right \|=x^{\beta-\mu-1}E_{\alpha, \beta-\mu}[J_{r_k}(\lambda_k)z],$$
where $\tilde e_n\equiv e_n^{\alpha,\beta-\mu}(\lambda,z)=z^n E_{\alpha, \alpha n+\beta-\mu}^{n+1}(\lambda z).$
Taking into account (\ref{eq191}), we get (\ref{eq52}).

3. From (\ref{eq51}) and (\ref{eq52}) it follows
\begin{equation} \label{eq53}
\left(D_{0x}^{\alpha}-A\right)x^{\beta-1}E_{\alpha, \beta}(Ax^{\alpha})=\frac{x^{\beta-\alpha-1}}{\Gamma(\beta-\alpha)}I.
\end{equation}


4. Let us act on the function $x^{\beta-1}E_{\alpha, \beta}(Ax^{\alpha})$ by the Dzhrbashyan -- Nersesyan operator.

We denote $\mu_j=\sum\limits_{i=0}^{j}\gamma_i.$
Then, by virtue of the formula (\ref{eq52}),
for $ \beta-\mu_j \geq 0$ we get the following formula
\begin{equation} \label{eq54}
D_{0x}^{\{\gamma_0, \gamma_1,...,\gamma_j\}}x^{\beta-1}E_{\alpha, \beta}(Ax^{\alpha})=
x^{\beta-\mu_j}E_{\alpha, \beta-\mu_j+1}(Ax^{\alpha}).
\end{equation}
In particular, for $\beta-\mu_j=0$ we have
\begin{equation} \label{eq55}
D_{0x}^{\{\gamma_0, \gamma_1,...,\gamma_j\}}x^{\beta-1}E_{\alpha, \beta}(Ax^{\alpha})=
E_{\alpha, 1}(Ax^{\alpha}).
\end{equation}
By using the formula (\ref{eq51}) from (\ref{eq54})  we obtain
$$D_{0x}^{\{\gamma_0, \gamma_1,...,\gamma_j\}}x^{\beta-1}E_{\alpha, \beta}(Ax^{\alpha})=
\frac{x^{\beta-\mu_j}}{\Gamma(\beta-\mu_j+1)}I+Ax^{\beta+\alpha-\mu_j}E_{\alpha,\beta+\alpha-\mu_j+1}(Ax^{\alpha})$$
for $\beta-\mu_j \geq 0.$

From the last equality, for $j=m, $ and $\gamma_i=\alpha_i$ $ (i = \overline{0, m}), $ we get
\begin{equation} \label{eq551}
\left(D_{0x}^{\{\alpha_0, ...,\alpha_m\}}-A\right)x^{\beta-1}E_{\alpha, \beta}(Ax^{\alpha})=
\frac{x^{\beta-\alpha-1}}{\Gamma(\beta-\alpha)}I,
\quad \beta\geq \mu_m=\alpha+1,
\end{equation}
and for $j=m,$ and $\gamma_i=\alpha_{m-i}$ $(i=\overline{0,m}),$ we get
\begin{equation} \label{eq56}
\left(D_{0x}^{\{\alpha_m, ...,\alpha_0\}}-A\right)x^{\alpha}E_{\alpha, \alpha+1}(Ax^{\alpha})=I.
\end{equation}

5. By virtue of formula (\ref{eq155}) we obtain
$$\lim\limits_{z \to 0} e_{n}^{\alpha, \beta}(\lambda, z)=
\left \{
\begin{array}{lr}
0,						 	& n>0,\\
 \frac{1}{\Gamma(\beta)},     	& n=0.
\end{array}
\right .
$$
From the last formula and (\ref{eq191}) we get
$$\lim\limits_{x\to 0}E_{\alpha, \beta}(A x^{\alpha})=  \frac{1}{\Gamma(\beta)}I,$$ 
and
\begin{equation} \label{eq555}
\lim\limits_{x\to 0}x^{\beta-1}E_{\alpha, \beta}(Ax^{\alpha})=
\left \{
\begin{array}{lr}
0,		&  \beta>1,\\
I,     	&  \beta=1.
\end{array}
\right .
\end{equation}

\begin{center}
{\bf 4. Green's formula for the operator $L$}
\end{center}

{\bf Lemma 1.}
{\it Let $f(x)\in L(0,l),$ then every regular solution to Problem 1 can be represented by formula (\ref{eq111}).
}

{\bf Proof.}
Let the function $V(x,t)$ be such that $D_{xt}^{\{\alpha_m, \alpha_{m-1},..., \alpha_{k}\}}V(x,t)\in C^1(0,l)\cup C[0,l]$  $(k=\overline{1,m})$ and 
 $D_{xt}^{\{\alpha_m, \alpha_{m-1},..., \alpha_0\}}V(x,t)\in C(0,l)\cup L(0,l)$ for any fixed $x\in [0,l].$ 
Using the definitions of the Dzhrbashyan -- Nersesyan and Riemann -- Liouville operators, the integration by parts formula, and the formula (\ref{eq11}), 
we transform the following integral
$$\int\limits_0^xV(x,t)D_{0t}^{\{\alpha_0, \alpha_1,...,\alpha_m\}}u(t)dt=
\int\limits_0^xV(x,t)D_{0t}^{\alpha_m-1}\frac{d}{dt}D_{0t}^{\{\alpha_0, \alpha_1,...,\alpha_{m-1}\}}u(t)dt=$$
$$=\int\limits_0^xD_{xt}^{\alpha_m-1}V(x,t)\cdot \frac{d}{dt}D_{0t}^{\{\alpha_0, \alpha_1,...,\alpha_{m-1}\}}u(t)dt=
S_1(x)-\int\limits_0^x\frac{d}{dt}D_{xt}^{\alpha_m-1}V(x,t)\cdot D_{0t}^{\{\alpha_0, \alpha_1,...,\alpha_{m-1}\}}u(t)dt=$$
$$=S_1(x)+\int\limits_0^xD_{xt}^{\alpha_m}V(x,t)\cdot D_{0t}^{\alpha_{m-1}-1}\frac{d}{dt}D_{0t}^{\{\alpha_0, \alpha_1,...,\alpha_{m-2}\}}u(t)dt=$$
$$=S_1(x)+\int\limits_0^xD_{xt}^{\alpha_{m-1}-1}D_{xt}^{\alpha_m}V(x,t)\cdot \frac{d}{dt}D_{0t}^{\{\alpha_0, \alpha_1,...,\alpha_{m-2}\}}u(t)dt=$$
$$=S_1(x)+\int\limits_0^xD_{xt}^{\{\alpha_m, \alpha_{m-1}\}}V(x,t)\cdot \frac{d}{dt}D_{0t}^{\{\alpha_0, \alpha_1,...,\alpha_{m-2}\}}u(t)dt,$$
where
$$S_1(x)=D_{xt}^{\{\alpha_m\}}V(x,t)\cdot D_{0t}^{\{\alpha_0, \alpha_1,..., \alpha_{m-1}\}}u(t)\big|_{t=0}^{t=x}.$$
Continuing similarly, we obtain
$$\int\limits_0^xV(x,t)D_{0t}^{\{\alpha_0, \alpha_1,...,\alpha_m\}}u(t)dt=
\sum\limits_{i=1}^{k}D_{xt}^{\{\alpha_m, \alpha_{m-1},..., \alpha_{m+1-i}\}}V(x,t)\cdot D_{0t}^{\{\alpha_0, \alpha_1,..., \alpha_{m-i}\}}u(t)\big|_{t=0}^{t=x}+$$
$$+\int\limits_0^xD_{xt}^{\{\alpha_m, \alpha_{m-1},..., \alpha_{m-k}\}}V(x,t)\cdot \frac{d}{dt}D_{0t}^{\{\alpha_0, \alpha_1,..., \alpha_{m-k-1}\}}u(t)dt=...=$$
$$=\sum\limits_{i=1}^{m}D_{xt}^{\{\alpha_m, \alpha_{m-1},..., \alpha_{m+1-i}\}}V(x,t)\cdot D_{0t}^{\{\alpha_0, \alpha_1,..., \alpha_{m-i}\}}u(t)\big|_{t=0}^{t=x}
+\int\limits_0^xD_{xt}^{\{\alpha_m, \alpha_{m-1},..., \alpha_0\}}V(x,t)\cdot u(t)dt.$$

Thus, the formula
$$\int\limits_0^x\left[V(x,t)D_{0t}^{\{\alpha_0, \alpha_1,...,\alpha_m\}}u(t)-D_{xt}^{\{\alpha_m, \alpha_{m-1},..., \alpha_0\}}V(x,t)\cdot u(t)\right]dt=$$
$$=\sum\limits_{i=1}^{m}D_{xt}^{\{\alpha_m, \alpha_{m-1},..., \alpha_{m+1-i}\}}V(x,t)\cdot D_{0t}^{\{\alpha_0, \alpha_1,..., \alpha_{m-i}\}}u(t)\big|_{t=0}^{t=x}$$
 holds.

Let the function $V(x,t)$ is a solution of the equation
\begin{equation} \label{eq31}
L^*V(x,t)\equiv D_{xt}^{\{\alpha_m, \alpha_{m-1},..., \alpha_0\}}V(x,t)-V(x,t)A=I,
\end{equation}
and satisfies following conditions
\begin{equation} \label{eq32}
\lim\limits_{t \rightarrow x} D_{xt}^{\{\alpha_m, \alpha_{m-1},..., \alpha_{m+1-i}\}}V(x,t)=0,\quad
 1\leq i\leq m.
\end{equation}
From the formulas (\ref{eq56}),  (\ref{eq555}) and
\begin{equation} \label{eq30}
D_{xt}^{\gamma}v(x-t)=D_{0y}^{\gamma}v(y)\big|_{y=x-t}
\end{equation}
it follows that the function
$$V(x,t)=(x-t)^{\alpha}E_{\alpha,\alpha+1}(A(x-t)^{\alpha})$$
is the solution of problem (\ref{eq31}), (\ref{eq32}).

From (\ref{sn}) and (\ref{eq31}) follows that
\begin{equation} \label{eq33}
V(x,t)Lu(t)-L^*V(x,t)\cdot u(t)=V(x,t)f(t)-u(t).
\end{equation}
After integrating the equality (\ref{eq33}), we obtain
\begin{equation} \label{eq34}
\sum\limits_{i=1}^{m}D_{xt}^{\{\alpha_m, \alpha_{m-1},..., \alpha_{m+1-i}\}}V(x,t)\cdot D_{0t}^{\{\alpha_0, \alpha_1,..., \alpha_{m-i}\}}u(t)\big|_{t=0}^{t=x}=
\int\limits_0^xV(x,t)f(t)dt-\int\limits_0^xu(t)dt.
\end{equation}
Using (\ref{eq11}) and (\ref{eq32}) from equality (\ref{eq34}) we get
\begin{equation} \label{eq35}
\int\limits_0^xu(t)dt=
\int\limits_0^xV(x,t)f(t)dt+
\sum\limits_{i=1}^{m}D_{xt}^{\{\alpha_m, \alpha_{m-1},..., \alpha_{m+1-i}\}}V(x,t)\big|_{t=0}u_0^{m-i}.
\end{equation}
Differentiating (\ref{eq35}), taking into account the equality $ V(x,x) = 0, $ we obtain
\begin{equation} \label{eq36}
u(x)=
\int\limits_0^xG(x,t)f(t)dt+
\sum\limits_{j=0}^{m-1}D_{xt}^{\{\alpha_m, \alpha_{m-1},..., \alpha_{j+1}\}}G(x,t)\big|_{t=0}u_0^j,
\end{equation}
where $G(x,t)=V_x(x,t)=G(x-t),$
$$G(x)=x^{\alpha-1}E_{\alpha, \alpha}\left(Ax^{\alpha}\right).$$
Applying to (\ref{eq36}) the formula
$$
\lim\limits_{t \rightarrow 0} D_{xt}^{\gamma}v(x-t)=D_{0x}^{\gamma}v(x),
$$
that follows from (\ref{eq30}), we get the equality (\ref{eq111}).
Lemma 1 is proved.

\begin{center}
{\bf 5. Proof of Theorem 1}
\end{center}

Let us prove that the function (\ref{eq111}) is the solution of Problem 1.
We denote the last term in right hand side of (\ref{eq111}) as $u_f(x)$ and
$u_C(x)=u(x)-u_f(x).$

By virtue of (\ref{eq54}) we get
\begin{equation} \label{eq561}
D_{0x}^{\{\alpha_m, ...,\alpha_k\}}G(x)=x^{\alpha-\nu_k}E_{\alpha, \alpha-\nu_k+1}(Ax^{\alpha})=x^{\mu_{k-1}-1}E_{\alpha, \mu_{k-1}}(Ax^{\alpha}),
\end{equation}
where $\mu_k=\sum\limits_{i=0}^{k}\alpha_i,$ $\nu_k=\sum\limits_{i=k}^{m}\alpha_i.$

Using the following equalities
$$\alpha-\mu_k-\nu_s+1=
\left \{
\begin{array}{ll}
\alpha_{s+1}+...+\alpha_{k-1}, 	& s<k-1, \\
 0, 						   		& s=k-1, \\
-(\alpha_k+...+\alpha_s), 			& s>k-1,
\end{array}
\right.$$
the formulas (\ref{eq51}), (\ref{eq561}) and (\ref{eq54}) we obtain
\begin{equation} \label{eq57}
D_{0x}^{\{\alpha_0, ...,\alpha_s\}}D_{0x}^{\{\alpha_m, ...,\alpha_k\}}G(x)=
\left \{
\begin{array}{ll}
x^{\alpha_{s+1}+...+\alpha_{k-1}}E_{\alpha, \alpha_{s+1}+...+\alpha_{k-1}+1}(Ax^{\alpha}), 	& s<k-1, \\
 E_{\alpha, 1}(Ax^{\alpha}), 						   		& s=k-1, \\
Ax^{\mu_{k-1}+\nu_{s+1}-1}E_{\alpha, \mu_{k-1}+\nu_{s+1}}(Ax^{\alpha}), 	& s>k-1.
\end{array}
\right.
\end{equation}

From (\ref{eq57}), (\ref{eq555}) and the inequality $\alpha_0 + \alpha_m> 1$ it follows that
\begin{equation} \label{eq58}
\lim\limits_{x\to 0}D_{0x}^{\{\alpha_0, ...,\alpha_s\}}D_{0x}^{\{\alpha_m, ...,\alpha_k\}}G(x)=
\left \{
\begin{array}{ll}
0, 		& s<k-1, \\
I,   		& s=k-1, \\
0, 	& s>k-1.
\end{array}
\right.
\end{equation}
Formula (\ref{eq58}) gives the relations
\begin{equation} \label{eq581}
\lim\limits_{x\to 0}D_{0x}^{\{\alpha_m, ...,\alpha_k\}}u_C(x)=u_{0}^{k}, \quad k=\overline{0, m-1}.
\end{equation}
For $ s = m $ from (\ref{eq56}) and (\ref{eq57}) we get
\begin{equation} \label{eq59}
\Bigl( D_{0x}^{\{\alpha_0, ...,\alpha_m\}}-A \Bigr)D_{0x}^{\{\alpha_m, ...,\alpha_k\}}G(x)=0.
\end{equation}

Thus, it remains for us to show that the last term in (\ref{eq111}) is a solution to Problem 1 with homogeneous conditions.

By virtue of the condition  $f(x)=D_{0x}^{\alpha_m-1}f_0(x),$ $f_0(x)\in L(0,l)$ and formulas
(\ref{eq11}) and (\ref{eq52}), we can write
$$u_f(x)=\int\limits_0^xG(x-t)D_{0x}^{\alpha_m-1}f_0(t)dt=
\int\limits_0^x G_0(x-t)f_0(t)dt,$$
where 
$$G_0(x)=D_{0x}^{\alpha_m-1}G(x)=x^{\mu_{m-1}-1}E_{\alpha,\mu_{m-1}}(Ax^{\alpha}).$$

The last equality yield
$$D_{0x}^{\{\alpha_0\}}u_f(x)=
\int\limits_0^x D_{xt}^{\alpha_0-1}G_0(x-t)\cdot f_0(t)dt.$$

Using formula (\ref{eq54}) we get
\begin{equation} \label{eq562}
D_{0x}^{\{\alpha_0, ...,\alpha_k\}}G_0(x)=x^{\mu_{m-1}-\mu_k}E_{\alpha, \mu_{m-1}-\mu_k+1}(Ax^{\alpha}).
\end{equation}
 
By virtue of (\ref{eq13}), (\ref{eq52}) and (\ref{eq555})  we have
$$D_{0x}^{\alpha_0}u_f(x)=\frac{d}{dx}D_{0x}^{\{\alpha_0\}}u_f(x)=
\int\limits_0^x D_{xt}^{\alpha_0}G_0(x-t)\cdot f_0(t)dt+$$
$$+\Bigl(\lim\limits_{t\to x}D_{0x}^{\alpha_0-1}G_0(x-t)\Bigr)f_0(x)=\int\limits_0^x D_{xt}^{\alpha_0}G_0(x-t)\cdot f_0(t)dt.$$

Continuing in a similar way, we obtain the following equalities
\begin{equation} \label{eq61}
D_{0x}^{\{\alpha_0,...,\alpha_k\}}u_f(x)=
\int\limits_0^xD_{xt}^{\{\alpha_0,...,\alpha_k\}}G_0(x-t)f_0(t)dt,
\end{equation}
and
$$\frac{d}{dx}D_{0x}^{\{\alpha_0,...,\alpha_k\}}u_f(x)=
\int\limits_0^x\frac{d}{dx}D_{0x}^{\{\alpha_0,...,\alpha_k\}}G_0(x-t)f_0(t)dt+$$
$$+\Bigl(\lim\limits_{t\to x}D_{xt}^{\{\alpha_0,...,\alpha_k\}}G_0(x-t)\Bigr)f_0(x)=
\int\limits_0^x\frac{d}{dx}D_{0x}^{\{\alpha_0,...,\alpha_k\}}G_0(x-t)f_0(t)dt.$$
For $k=m-1$ from (\ref{eq61}) we get
\begin{equation} \label{eq62}
D_{0x}^{\{\alpha_0,...,\alpha_{m-1}\}}u_f(x)=\int\limits_0^x E_{\alpha,1}(A(x-t)^{\alpha})f_0(t)dt.
\end{equation}
From the relation (\ref{eq62}) using (\ref{eq51}) we get
$$\frac{d}{dx}D_{0x}^{\{\alpha_0,...,\alpha_{m-1}\}}u_f(x)=
\frac{d}{dx}\int\limits_0^x (x-t)^{\alpha}E_{\alpha,\alpha+1}(A(x-t)^{\alpha})f_0(t)dt=$$
$$=\lim\limits_{t\to x}E_{\alpha,1}(A(x-t)^{\alpha})f_0(t)+\int\limits_0^x (x-t)^{\alpha-1}E_{\alpha,\alpha}(A(x-t)^{\alpha})f_0(t)dt=$$
$$=f_0(x)+A\int\limits_0^x G(x-t)f_0(t)dt.$$
From the last relation  follows that
$$D_{0x}^{\{\alpha_0, ...,\alpha_m\}}u_f(x)=D_{0x}^{\{\alpha_m\}}\frac{d}{dx}D_{0x}^{\{\alpha_0, ...,\alpha_m-1\}}u_f(x)=$$
$$=D_{0x}^{\alpha_m-1}\Bigl( f_0(x)+A\int\limits_0^x(x-t)^{\alpha-1}E_{\alpha,\alpha}(A(x-t)^{\alpha})f_0(t)dt\Bigr)=$$
\begin{equation} \label{eq63}
=f(x)+A\int\limits_0^xG(x-t)f(t)dt=f(x)+Au_f(x).
\end{equation}
Equalities (\ref{eq562}) and (\ref{eq61}) lead us to the relation
\begin{equation} \label{eq64}
\lim\limits_{x\to 0}D_{0x}^{\{\alpha_0,...,\alpha_k\}}u_f(x)=0,
\quad 0\leq k\leq m-1.
\end{equation}

The relations (\ref{eq581}), (\ref{eq59}), (\ref{eq63}) and (\ref{eq64}) mean
that the function (\ref{eq111}) is the solution to Problem 1.
The uniqueness of the solution to Problem 1 follows from Lemma 1.

Theorem 1 is proved.

\begin{center}
{\bf 5. Conclusions}
\end{center}

The article investigates the initial problem for a linear system of ordinary differential equations with the Djrbashyan -- Nersesyan fractional differentiation operator with constant coefficients.
To solve the problem under study, the Green's function method is implemented, an integral representation of the solution is obtained,
the properties of the matrix Mittag-Leffler function (autotransformation formula, fractional integro-differentiation formulas, etc.) are investigated.
The results obtained can be used to study local and nonlocal boundary value problems for system (1).

Conflicts of Interest: The author declares no conflict of interest in this paper.


\vskip 3mm

{\it Submitted July 27, 2020}

\vskip 3mm

{{\footnotesize 
Mamchuev Murat Ocmanovich

Institute of Applied Mathematics and Automation KBSC RAS

89~A Shortanov Street, Nalchik, 360000, Russia

mamchuev@rambler.ru
}}

\end{document}